\documentclass[11pt]{amsart}
\usepackage{amsfonts, eucal, amssymb}
\usepackage[all]{xy}
\begin{document}

\newcommand{\E}{{\mathbb E}}
\newcommand{\h}{{\hbar}}
\newcommand{\Oh}{{\mathcal O}}
\newcommand{\oh}{{\bf o}}
\newcommand{\Mor}{{\rm Mor}}       
\newcommand{\Lie}{{\rm Lie \;}}
\newcommand{\Maps}{{\rm Maps}}
\newcommand{\Trans}{{\rm Trans}}
\newcommand{\Z}{{\mathbb Z}}
\newcommand{\Aut}{{\rm Aut}}
\newcommand{\bPhi}{{\boldsymbol{\Phi}}}
\newcommand{\Fix}{{\rm Fix}}
\newcommand{\Iso}{{\rm Iso}}
\newcommand{\N}{{\mathbb N}}
\newcommand{\cN}{{\mathcal N}}
\newcommand{\C}{{\mathbb C}}
\newcommand{\A}{{\mathcal A}}
\newcommand{\G}{{\mathcal G}}
\newcommand{\mG}{{\mathbb G}}
\newcommand{\D}{{\mathcal D}}
\newcommand{\Diff}{{\rm Diff}}
\newcommand{\Sets}{{\rm Sets}}
\newcommand{\Spaces}{{\rm Spaces}}
\newcommand{\Vect}{{\rm Vect}}
\newcommand{\dis}{{\rm discrete}}

\title{Theories of anything}
\author{Jack Morava}
\address{The Johns Hopkins University,
Baltimore, Maryland 21218}
\email{jack@math.jhu.edu}
\thanks{This work was supported by the NSF}
\date{20 July 2012}
\begin{abstract}{We suggest a generalization of $\pi_0$ for topological groupoids, which
encodes incidence relations among the strata of the associated quotient object, and argue for its
utility by example, starting from the orbit categories of the theory of compact Lie groups.
\medskip

\noindent
One of the points of this note is that Thom's theory of structurally stable forms fits
quite nicely with the categorical theory of databases developed recently by D. Spivak; the other is
that the stratifications studied by Thom are closely related to the phase transitions studied in
physics, and that the generalization of $\pi_0$ proposed here may be useful in their study: in
particular, in organizing our understanding of the scaling laws which naturally accompany such
phenomena, in the theory of condensed matter, biology, finance \dots.}\end{abstract}

\maketitle

\section{Introduction}\bigskip

\noindent
Many difficulties in mathematical classification involve degenerate objects: the theory of conic
sections is a classical example, and geometric invariant theory has evolved powerful techniques
to deal with such questions. A central issue, in topological terms, is that an action of a reasonable
group $G$ on a nice space $X$ can quite easily fail to have a Hausdorff quotient space. \bigskip

\noindent
This paper is a rough and intuitive sketch. It  proposes encoding the incidence structure of such
bad quotients in a database category $\bPhi_0$ which generalizes the classical notion of the set
$\pi_0$ of components of a topological space. It presents some examples, and argues for the
hope of a coherent framework to accomodate them; it is thus mostly about what physicists call
phenomenology. \bigskip

\noindent
We focus on {\bf smooth transformation groupoids} $[X/G]$, perhaps not finite-dimensional, though 
the stratified spaces of Thom [29] (and more general topological groupoids [12, 23, 30]) share many 
features with this class of objects. An archetypical example is Arnol'd's classification
of isolated singularities of holomorphic functions; the bulk of this paper [\S 3] is concerned with
that case, which illustrates particularly clearly that (in physical examples) degeneracy relations
among orbits can lead to the existence of scaling laws near phase transitions. \bigskip

\noindent
{\bf Acknowledgements:} This work began in conversations with David Spivak; I will try to
explain below why the `phase diagrams' defined here fit into his thinking [26, 27] about
databases. I also owe thanks to Jos\'e Manuel Gomez and John Lind, for patient conversations
about equivariant topology, and, at a different coherence length, to Ben Mann, Simon Levin,
Andrew Salch, Jim Stasheff, and Abel Wolman III. \bigskip

\section{Transformation groupoids and orbit categories} \bigskip

\noindent
{\bf 2.1} If $G$ is a compact Lie group, there is a topological category whose objects are its
closed subgroups, eg $ H_0,H_1,\dots$ , with the spaces
\[
G-\Maps(G/H_0,G/H_1) \; := \: \Mor(H_0,H_1)       
\]
of equivariant maps as morphism objects. It is {\bf equivalent} to the topological category
$\Oh(G)$ with the (totally disconnected [10, 13]) space of conjugacy classes of closed
subgroups as its space of objects, and
\[
\Mor_\Oh(H_0,H_1) \; := \; \Trans_G(H_0,H_1)/H_0
\]
as morphisms; where 
\[
\Trans_G(H_0,H_1) \; = \; \{g \in G \:|\: gH_0g^{-1} \subset H_1 \}
\]
is the `transporter' of $H_0$ to $H_1$ in $G$. Thus 
\[
\Aut_\Oh(H) \;  = \; N_G(H)/H \; := \: W_G(H)
\]
is a kind of Weyl group. The space of subgroups of $G$ is filtered by dimension, and the space
of morphisms from $H_0$ to $H_1$ will be empty if dim $H_0 >$  dim $H_1$. \bigskip

\noindent
Many of the constructions of this paper have analogs for locally compact groups, using the
Chabauty topology [5]. \bigskip                        

\noindent
{\bf Definition} $\Oh_0(G)$ is the (topological) category with conjugacy classes of
closed subgroups of $G$ as objects, and the discrete spaces
\[
\Mor_{\Oh_0}(H_0,H_1) \; := \; \pi_0 \Mor(H_0,H_1) 
\]
as morphism objects. \bigskip

\noindent
{\bf 2.2} The topological category $[X/G]$ defined by an action
\[
G \times X \to X
\]
of $G$ on a space $X$ has elements $x_0,x_1,\dots \in X$ as objects, with
\[
\Mor_{[X/G]}(x_0,x_1) \; := \; \{g \in G \:|\: gx_0 = x_1\}
\]
as morphisms; in particular, if $X$ is at all nice, the isotropy group
\[
\Aut_{[X/G]}(x) \; = \; \{ g \in G \:| \: gx = x \} \; := \; \Iso(x)
\]
of $x$ is closed, and
\[
\Fix(X) : H \mapsto \{x \in X \:| hx = x \; \forall h \in H \} \; 
\]
defines a contravariant functor from closed subgroups of $G$ to spaces. Taking components
defines a presheaf
\[
H \mapsto \pi_0 \Fix(H) : \Oh_0(G) \to (\Sets) \;.
\] 

\noindent
{\bf  Definition:} The pullback (or coend [19]) 
\[
\xymatrix{\bPhi_0[X/G] \ar@{.>}[r] \ar@{.>}[d] & (\Sets)_* \ar[d] \\
\pi_0\Oh(G) \ar[r]^{\pi_0 \Fix(X)} & (\Sets) }
\]
defines Grothendieck's category of elements for the functor $\pi_0 \circ \Fix(X)$. Its objects are 
pairs $(c,H)$, consisting of (conjugacy classes of) closed subgroups $H$, together with a
component $c$ of $\Fix(H)$; thus 
\[
\Aut_\bPhi(c,H) \; = \; \{ g \in \pi_0 W_G(H) \:|\: gc = c \} \;.
\]
This defines a functor $\bPhi_0$ from transformation groupoids to (Boolean) topological
categories, which I'd like to call the {\bf phase diagram} [25] of the groupoid. \bigskip

\noindent
{\bf Example} $\bPhi_0[1/G]  =  \Oh_0(G)$ is the discretized orbit category defined
above, and the forgetful functor
\[
\bPhi_0[X/G] \to \bPhi_0[1/G]
\]
presents its domain as a category fibered over its target. \bigskip

\noindent
{\bf 2.3} In fact $[X/G]$ can be expressed as the (topological) category of elements of the 
$G$-space $X$, regarded as a functor from $[1/G]$ to topological spaces. It can similarly be 
identified with the category of elements of the presheaf
\[
\Fix(X) : \Oh(G) \to (\Spaces) \;;
\]
which implies the existence of a functor
\[
[X/G] \to \bPhi_0[X/G]
\]
generalizing the classical map $X \to \pi_0 X$. \bigskip

\noindent
Classically, $\pi_0$ is adjoint
\[
\Maps_\Spaces(X,S_\dis) \; = \; \Maps_\Sets(\pi_0X,S)
\]
to the inclusion of sets into spaces, and it would be very useful to characterize $\bPhi_0$ in
similar terms. Peter May observes [p.c., but cf [20]] that the category of finite topological
spaces is Cartesian closed; it is a little small for our purposes, but one can ask if the
the inclusion of the category of profinite topological categories into some reasonable larger class
of topological categories might possess a (2-)adjoint. \bigskip

\noindent
{\bf 2.4} When $X$ is a manifold, and $G$ acts smoothly on it, I'll refer to $[X/G]$ as a
smooth transformation groupoid. In the finite-dimensional case, the {\bf slice theorem} [14, 18]
implies that any $x \in X$ has a $G$-invariant open neighborhood $S(x)$ equivariantly
diffeomorphic to
\[
\cN_x \times_{\Iso(x)} G \;,
\]
where $\cN_x$ is essentially [2 I \S 1.7] the linear representation of $\Iso(x)$ defined by the
cokernel of the tangent map to the inclusion of the $G$-orbit through $x$ into the total space
$X$. A morphism
\[
\alpha \in \Mor_\bPhi((x_0,H_0),(x_1,H_1))
\]
defines (the germ of) an embedding
\[
S(x_0) \subset S(x_1)
\]
and hence (the conjugacy class of) a monomorphism
\[
\cN_{x_0} \to \cN_{x_1} \;,
\]
equivariant with respect to the embedding $H_0 \to H_1$ defined by the adjacency $\alpha$.
\bigskip

\noindent
{\bf Definition} I'll call the resulting functor
\[
\cN : \bPhi_0[X/G] \to \Vect
\]
(with values in the category of real vector spaces and monomorphisms), the {\bf degeneracy
quiver} of $[X/G]$. \bigskip

\section{Some motivating examples} \bigskip

\noindent
{\bf 3.1 Arnol'd's category} of (stable equivalence classes of) isolated singularities of
holomorphic functions: \bigskip

\noindent
If $f : (\C^n,0) \to (\C,0)$ is the germ of a holomorphic function, its associated algebra $\oh_f$ is
the quotient of the local ring $\C\{x_1,\dots,x_n\}$ of functions holomorphic at 0 by the
(Jacobian) ideal of $f$ generated by its partial derivatives $\partial f/\partial x_i$. \bigskip

\noindent
If the Milnor number $\mu(f) := \dim_\C \oh_f $ is finite, $f$ defines an {\bf isolated} singularity 
at 0 [2 I \S 1.4]. Let $\A_n$ denote the space of such function-germs, and let $\G_n$ be the
group (under composition) of germs of biholomorphic maps $(\C^n,0) \to (\C^n,0)$, which acts
on $\A_n$ by left composition. \bigskip

\noindent
The transformation groupoid $[\A_n/\G_n]$ is smooth but infinite-dimensional;  but since the
Milnor numbers of its objects are finite, their orbits are determined by the behavior of some
representative jet [2 I \S 1.5]. The resulting groupoid therefore looks in many ways like a global
quotient defined by an action of a locally compact group. The stabilization map
\[
f \mapsto f + x_{n+1}^2 : \C\{x_1,\dots,x_n\} \to \C\{x_1,\dots,x_{n+1}\}
\]
[2 I \S 1.3] defines Arnol'd's big stack
\[
\lim_{n \to \infty} [\A_n/\G_n] \; := \; [\A/\G]
\]
of (stable equivalence classes of) isolated singularities of holomorphic functions. It is naturally
stratified by orbit codimension; the {\bf modality}
\[
m(f) \; := \: (\mu(f) - 1) - {\rm codim} \; f
\]
is an important related (subtle) invariant. \bigskip

\noindent
Arnol'd and his school have worked out the structure of this stack, for $m < 3$; their results seem
to me a plausible candidate for an {\bf olog} in the sense of [27], ie a kind of annotated database.
In particular, the stratification of $[\A/\G]$ up to codimension ten is well-understood; see for
example the adjacency diagrams in [2 I \S 2]. \bigskip

\noindent
{\bf 3.2} The moduli stack of {\bf one-dimensional formal groups} is how I got into this [22 \S
1.4]; it is a big arithmetic stack (ie defined over Spec $\Z$), with a nontrivial but very simple
stratification (isomorphic to $\N$, regarded as a poset with its usual order) at each prime. Formal
groups of higher dimension are of course much more complicated \dots \bigskip

\noindent
{\bf 3.3 Knots} (and links \dots) can be studied similarly, as topological groupoids defined by a
space 
\[
X \; := \; {\rm Imm}(S^1,S^3)
\]
of immersions, under the action of the group 
\[
G \; := \; \Diff(S^1) \times \Diff(S^3)
\]
of diffeomorphisms of the domain and range; the theory of finite-type invariants comes from its
stratification by self-intersection number. Knot tables (cf eg [3, 29]) are another natural class 
of examples of mathematical ologs. \bigskip

\noindent
{\bf 3.4.1} Similarly, {\bf Riemann surfaces} seem to fit well enough into some such framework;
but formulating the moduli problem as a global quotient is not so easy [15], so I'll leave this
example aside; however, \bigskip

\noindent
{\bf 3.4.2} The topological groupoid of Riemannian metrics on a manifold, under the action of
its group of diffeomorphisms, is the natural configuration space for general relativity [9] (just
as the analogous stacks of connections on principal bundles are relevant to gauge theory). Classical 
solutions of the Einstein equations tend to have large symmetry groups, and thus represent quite singular
points of this moduli stack \dots \bigskip

\noindent
{\bf 3.5} Finally, Thom's notion [30 \S 2.1]  of a {\bf structurally stable form} lies behind this
whole essay. A cartoon version of his program asks for the structure of the big stack of
everything in the world, modulo the pseudogroup defined by their local isomorphisms. The
related theory of stratified sets, among other things, is one of his creations (see further 
\S 5.2), but in his day the language of (higher) category theory was only beginning to emerge. 
One of the points of this note is the close compatibility of that language (and with Spivak's 
theory of databases) with Thom's thinking. \bigskip

\noindent
$\; \bullet$ A very natural, highly nontrivial example of a Thomist olog, dear to my heart, is the
stratification of the space of everything in the world defined by the Bantu noun-class system [6,
8]). The `Leitfaden' in Serre's {\bf Corps locaux} [p 13] is another striking example \dots
\bigskip

\section{Adjacency and phase transitions} \bigskip

\noindent
{\bf 4.1} Arnol'd's simple singularities (whose orbits have no moduli, ie $m = 0$, eg with
codimension $\leq 7$  [2 II \S 2.4]) are worth further discussion. [They extend Thom's
classification [30 \S 3.2, 5.2-5.4] of elementary catastrophes.] In this class of examples, the 
nerve of the groupoid $[\cN_f/\Iso(f)]$ is the classifying space of a generalized braid group 
[2 II \S 5.3] and the adjacency maps of \S 1.4 above can be described in terms of inclusions of 
Dynkin diagrams [2 II \S 5.9, 2 \S 6 (Fig 39)]. \bigskip

\noindent
Functions with simple singularities have nice quasihomogeneous normal forms, which endow
their associated algebras with canonical Euler derivations [2 II \S 5.7]. A degeneration $\alpha : f
\mapsto g$ (ie with adjacent orbits [2 I \S 2.7]) defines a homomorphism
\[
(\oh_f,\D_f ) \to (\oh_g,\D_g) 
\]
of differential algebras; its cokernel $\cN^f_g$ can be identified with that of the map defined
in \S 2.4. In particular, when the codimension jump from $f$ to $g$ is one, we can think of 
$\cN^f_g$ as generated by $\D_g$ (modulo $\D_f$).\bigskip

\noindent
{\bf 4.2.1} Physical systems are often analyzed in Morse-theoretic terms, by defining a suitable
Lagrangian (action) functional of on some space $X$ of states, perhaps invariant under some
(large) symmetry group $G$ (eg as in gauge theory). Such models can be formulated in terms of 
the gradient flow associated to the Lagrangian function on the associated quotient object; and in
interesting cases these quotients are not manifolds, but are instead stratified sets. \bigskip

\noindent
These stratifications can be described in terms of (`spontaneously') broken symmetry: that is, the
states of the system may have more (or less) symmetry than the equations of motion themselves
[30 \S 6.1 C,E]\begin{footnote}{In statistical mechanics, this goes back to Landau's work in
the late 1930's, where it is related to failure of the partition function to be analytic.}
\end{footnote}. In this conceptualization, phase transitions correspond to arrows
\[
(x_0,H_0) \to (x_1,H_1) \in \bPhi_0[X/G]
\]
which identify $H_0$ with a proper subgroup of $H_1$.  In particular, when dim $H_1 >$ dim
$H_0$, the normal bundle $\cN^{x_0}_{x_1}$ parametrizes the loss or gain in infinitesimal
symmetry associated to the transition $x_0 \to x_1$. \bigskip

\noindent
A classical theorem (perhaps better: a classical principle) of Noether associates a conserved
quantity to a continuous symmetry of a mechanical system; conservation of momentum, for
example, is a consequence of translation invariance, and angular momentum is similarly related
to rotational invariance. In the case above we expect to see, not a conserved quantity but rather
the appearance of an {\bf order parameter} $\eta$ [1 p 449, 16], which vanishes as $H_0
\to H_1$. In this sense $\cN^{x_0}_{x_1}$ is the Lie algebra generated by $\partial/\partial
\eta$ (when $\dim_\C \Lie H_1/\Lie H_0 = 1$). \bigskip

\noindent
{\bf 4.2.2} For adjacent simple singularities in $[\A/\G]$, the derivation $\partial/\partial 
\eta$ can be identified with the relative Euler operator above. In classical systems with 
$g \in \A$ as potential function, this generates an action of $\mG_m \cong \C^\times$: the {\bf 
renormalization group} [4, 17 \S 5.1] associated to the phase transition. \bigskip

\noindent
Eigenfunctions of Euler operators satisfy scaling laws, eg like those satisfied by `fat-tailed'
probability distributions [11, 17]; see further \S 5.4. In fact \bigskip

\begin{quotation}{\bf Rings of functions of observables defined near phase transitions are
naturally graded by scaling laws.}\end{quotation} \bigskip

\noindent
This is familiar in physics, but may have broader applications, eg to the complex systems
encountered in fields such as biology and finance. The second point of this paper is to suggest
that constructions such as $\bPhi_0$ may be useful in organizing our thinking about these scaling
laws. \bigskip 

\section{Afterthoughts} \bigskip

\noindent
{\bf 5.0} Jim Stasheff reminds me of Thom's dictum, that `Une th\'eorie qui explique tout 
n'explique rien'. The working title of these notes was `Phase diagrams for big stacks';
it is an attempt, not at a theory of anything, but of a theory of {\bf theories} of 
anything. \bigskip

\noindent
{\bf 5.1} (re \S 2.3) Substantial issues of finiteness, ultimately (point-set) topological [cf eg [12
\S 2.2]], lie behind the notion of a `big stack'; I have left this term undefined in the hope that
someone with a better understanding of such things will take up the challenge. \bigskip

\noindent
Examples in different contexts have been successfully treated by
various methods (eg in terms of pro-algebraic structures in ex's
3.1-2 above, while Sobolev `riggings' [9] are useful in examples
like 3.4.2) and it may be that no easy general formalism accomodates
them all; but the existence of slices [23] seems to be a very useful
general property. [A slice at an element $x$ of a topological groupoid
is presumably a nice neighborhood of its component in its equivalence
class; but `nice' will require further specification
in infinite-dimensional contexts.]\bigskip

\noindent
Finally, the 2-categorical invariance of finiteness conditions
seems quite subtle; this brings to mind G. Segal's remark, that
St. Anselm's ontological argument is basically Zorn's lemma
without chain conditions. \bigskip

\noindent
{\bf 5.2} (re \S 3.5) A stratified set $X$, eg an algebraic variety, is a union
of locally closed subsets $X_i$ indexed by a poset $I$, perhaps endowed with a monotone
increasing integer-valued codimension function. The category of elements of the functor
\[
i \mapsto \pi_0(\overline{X_i}) : I \to ({\rm Sets}) 
\]
defines a category $\bPhi_0(X)$ which (surprisingly) seems to have no name in algebraic 
geometry. \bigskip

\noindent
{\bf 5.3} (re \S 4.2.1) In statistical mechanics, the Helmholtz free energy of a system with 
Langrangian $L$ is also\begin{footnote}{More precisely, it's $\theta^{-1} \Gamma(\theta)$, with 
$\theta^{-1} = -kT$, where $T$ is the temperature and $k$ is Boltzmann's constant. In
population biology $\theta^{-1}$ is the number of organisms [7], and in finance \dots}\end{footnote} 
the cumulant generating function
\[
\Gamma(\theta) = \log \E (\exp(\theta L))
\]
of the random variable $L$. The stationary phase principle for the Feynman-Gibbs
partition function
\[
Z(\h) = \exp \Gamma(i\h^{-1}) = \int \exp(iL(x)/\h)dx
\]
suggests the interest of asymptotic expansions of 
\[
\Gamma(i\h^{-1}) = \log Z(\h) \in \C((\h)) \;,
\]
but in one of the simplest cases, ie the Bernoulli process with
\[
\Gamma(\theta) = \log (1 - p(1 - e^\theta)) \;,
\]
this is unproductive. On the other hand, in evolutionary contexts the Cram\'er `internal 
energy' (defined by the Legendre transform $\Gamma^*$ [21, 24] of $\Gamma$, which measures
the rate of the system's excursions from the mean [32]), may be more relevant. In this example
we have 
\[
C(x) := -\Gamma^*(x) = - x \log (1-p) + \log(1-p) + S(x) \;,
\]
with 
\[
S(x) = - x \log x - (1-x) \log (1-x) \;,
\]
and hence
\[
C(i\h^{-1}) \sim - i\h^{-1} \log(1-p) + \log \h + [\log(1-p)) - 1 + i \pi/4]  - \sum_{n\geq 1}
\frac{(i\h)^n}{n(n+1)} \; \dots
\]
[It is in some sense well-known that Planck's constant $h$ is {\bf not a number}. Geometric
quantization identifies $h^{-1}$ with the Chern class of a Hermitian line bundle
$(L,\nabla)$ with connection on the phase space $X$ of a physical system, ie
\[
i\h^{-1} \sim [\omega] \in H^2(X,2\pi i\Z) \subset H^2_{\rm dR}(X,\C) \;;
\]
from this point of view, $2 \pi i$ is not a number either: it is a motive\dots]
\bigskip

\noindent
{\bf 5.4} (re \S 4.2.2) It is therefore not surprising to find abundant non-Gaussian
phenomena in the vicinity of phase transitions; this may lie behind the biologists'
intuition, that living organisms colonize, and thrive, in systems on the `edge of
chaos'. It is worth noting the existence of a robust theory of non-Gaussian but 
`stable' probability distributions (going back to work of Kolmogorov and
Gnedenko in the 1930's [28]) which need not possess higher moments. \bigskip 
\bigskip

\bibliographystyle{amsplain}

\end{document}